\documentclass[11pt]{article}
\usepackage{amsmath}
\usepackage[all]{xy}
\usepackage{mathrsfs}
\usepackage{amssymb}
\usepackage{amsthm}
\usepackage[pdftex]{graphicx}
\numberwithin{equation}{section}

\begin{document}
\begin{center} {\bf {\LARGE Achievable Sets in $\mathbb{Z}^n$}}
\end{center}
\vskip 0.2in
{\small ABSTRACT: What sets $A \subset \mathbb{Z}^n$ can be written in the form $(K-K) \cap \mathbb{Z}^n$, where $K$ is a compact subset of $\mathbb{R}^n$ such that $K+\mathbb{Z}^n=\mathbb{R}^n$?  Such sets $A$ are called achievable, and it is known that if $A$ is achievable, then $\langle A \rangle=\mathbb{Z}^n$.  This condition completely characterizes achievable sets for $n=1$, but not much is known for $n \ge 2$.  We attempt to characterize achievable sets further by showing that with any finite, symmetric set $A \subset \mathbb{Z}^n$ containing zero, we may associate a graph $\mathcal{G}(A)$.  Then if $A$ is achievable, we show the set associated to some connected component of $\mathcal{G}(A)$ is achievable.  In two dimensions, we can strengthen this theorem further.  Further generalizations and open questions are discussed.  Throughout, the language and formalism of algebraic topology are useful.}

\vskip 0.4in

\begin{section}{Introduction}
\vskip 0.2in
Consider a subset $K$ of $\mathbb{R}^n$ with the property that $K$ is compact and $K+\mathbb{Z}^n=\mathbb{R}^n$.  Such a set is called an $\mathcal{N}$-\emph{set}.  In \cite{3}, Nathanson asked what subsets of $\mathbb{Z}^n$ can be written in the form $(K-K) \cap \mathbb{Z}^n$ for $K$ an $\mathcal{N}$-\emph{set}.  Such a subset is called \emph{achievable}.  Nathanson answered this question for $n=1$ by showing that $$A \subset \mathbb{Z} \hskip 0.1in \text{is \hskip 0.03in achievable} \iff \langle A \rangle=\mathbb{Z}.$$
The question has a natural generalization where $\mathbb{R}^n$ is replaced by a metric space $X$ and $\mathbb{Z}^n$ is replaced by a group $G$ which acts on $X$ so that the group action is \emph{geometric}. This notion is discussed in depth the appendix of \cite{3}.  Nathanson showed that even in the more general case (which includes the $n \ge 1$ case above), if $A \subseteq G$ is an achievable set, then $A$ generates $G$.
\vskip 0.2in
In the general case, this condition is not sufficient.  Ljujic and Sanabria \cite{2}, and Borisov and Jin \cite{1}, independently prove that any achievable set in $\mathbb{Z}^2$ contains elements with both coordinates zero, proving that the above condition is not sufficient even in two dimensions.  So then what additional conditions are necessary for a set to be achievable?  Can we find a nice characterization of achievable sets in $\mathbb{Z}^n$, even when $n=2$?  We focus on this question.
\vskip 0.3in
In Section 2, we discretize the problem, as is done in \cite{2} and \cite{1}.  In Section 3, we develop language and draw an analogy to the cohomology of the torus, so that in the next section, we may naturally prove Theorem 4.2, which is a statement about the `connectedness' of achievable sets.  In Section 5, we give an example of how this theorem gives us a good sense of what sets are achievable in $\mathbb{Z}^2$, and then strengthen the theorem in two dimensions (Theorem 5.3).
\vskip 0.2in
In Section 6, we state a generalization of the problem itself, and then in the context of this new problem, state Conjecture 6.3, a generalization of Theorem 4.3.  Theorem 4.3 gives us a better sense of achievable sets in two dimensions, but Conjecture 6.3 gives us a better sense of what sets are achievable in $n$ dimensions for all $n \ge 2$.  In Section 7, we conclude with further directions for research and open questions.
\end{section}
\vskip 0.5in
\begin{section}{A Discretization Step}
We begin by proving a general lemma which, in the case of $X=\mathbb{R}^n$, and $G=\mathbb{Z}^n$, allows us to discretize the problem.
\vskip 0.3in
{\bf Lemma 2.1.} \emph{Let $K_0 \supseteq K_1 \supseteq \cdots \supseteq K$ be a sequence of compact sets in $X$ such that $\bigcap\limits_{i=0}^{\infty}K_i=K$.  Then for some $m$, $A(K_m)=A(K)$.}
\vskip 0.1in
{\bf Proof:} Let $A_i$ denote $A(K_i)$ for $i=0, 1, \ldots$.  Because the group action is discontinuous, $A_i$ is finite for each $i$.  It is also clear that $A_0\supseteq A_1 \supseteq \cdots \supseteq A(K)$.  Let us suppose that each $A_i$ properly contains $A(K)$.  Then it follows that there is some $g$ which is contained in every $A_i$, but not contained in $A(K)$.  Hence, for $i=0, 1, \ldots$, we have $x_i, y_i \in K_i$ so that $g(x_i)=y_i$.
\vskip 0.2in
The sequence of pairs $(x_i, y_i)$ is contained in $K_0 \times K_0$, which is a compact subset of $X \times X$.  Since $X$ is a metric space, $X \times X$ is a metric space under the supremum norm.  Hence, $K_0 \times K_0$ is sequentially compact, and so the sequence $\{(x_i, y_i)\}_{i \ge 0}$ has a convergent subsequence $(x_{a_1}, y_{a_1}), (x_{a_2}, y_{a_2}), \ldots$.  Call its limit $(x, y)$.
\vskip 0.2in
Because this sequence converges under the supremum norm, we have that the sequence $x_{a_1}, x_{a_2}, \ldots$ converges to $x$, and the sequence $y_{a_1}, y_{a_2}, \ldots$ converges to $y$.  Notice that $x$ is contained in $K_{a_1}, K_{a_2}, \ldots$, so it must be contained in their intersection, $K$.  Similarly, $y \in K$.
\vskip 0.2in
Lastly, note that $d(x_i, x)=d(g(x_i), g(x))$, because the group action is isometric.  Hence, as $i$ goes to infinity, $y_{a_i}$ gets arbitrarily close to $g(x)$, and so $g(x)=y$.  Hence, it follows that $g \in A(K)$, which is a contradiction.  Hence, it follows that no such $g$ exists, and so $A_m=A(K)$ for some $m$ sufficiently large. $\blacksquare$
\vskip 0.4in
We now apply this lemma to a `nice' sequence of compact sets.
\vskip 0.3in
{\bf Definition 2.2.} Let $i_1, \ldots, i_n \in \mathbb{Z}$, $m \in \mathbb{N}$.  Consider the closed sets $$S_m(i_1, \ldots, i_n)=\left[\frac{i_1}{m}, \frac{i_1+1}{m}\right] \times \cdots \times \left[\frac{i_n}{m}, \frac{i_n+1}{m}\right]$$ in $\mathbb{R}^n$.  Call a set \emph{$m$-discrete} if it can be written as a finite union of sets $S_m(i_1, \ldots, i_n)$.  Call a set \emph{exactly $m$-discrete} if it can be written as a union of exactly $m^n$ of these hypercubes, and is an $\mathcal{N}$-\emph{set}.  Call a subset of $\mathbb{R}^n$ \emph{discrete} if it is $m$-discrete for some $m$, and correspondingly define an \emph{exactly discrete set}.
\vskip 0.3in
{\bf Proposition 2.2.} \emph{Let $A \subset \mathbb{Z}^k$ be achievable.  Then there is an exactly discrete $\mathcal{N}$-\emph{set} $K$ such that $A=A(K)$.}
\vskip 0.2in
\emph{Proof:} Let $(K-K) \cap \mathbb{Z}^n=A$ for some $K$, and for each $m$, let $K_m$ be the minimal covering of $K$ by sets $S_{2^m}(i_1, \ldots, i_n)$.  Then clearly $K_0 \supseteq K_1 \supseteq \cdots \supseteq K$.  Then by Lemma 2.1, there is some $m$ such that $(K_m-K_m) \cap \mathbb{Z}^n=A$.
\vskip 0.2in
Now notice that if a set $B$ is achievable by an exactly discrete set $K$, then for any $x$, $B \cup \{\pm x\}$ is achievable by an exactly discrete set $K'$.  Take the central $\frac{1}{3^n}$-th of one of the small hypercubes in $K$, and translate it by $x$.  This yields a $3m$-discrete set $K'$  such that $(K'-K') \cap \mathbb{Z}^n=B \cup \{\pm x\}$.  It follows that if $A$ is achievable by discrete set $\mathcal{N}$-\emph{set}, then it is achievable by an exactly discrete set as follows.  Throw away squares in this discrete set until we have an exactly discrete set, and then add back in any differences which were lost in this process by the above algorithm.  It follows that if $A$ is achievable, it is achievable by an exactly discrete set. $\blacksquare$
\vskip 0.2in
So now we may restrict our attention to these exactly discrete sets.  We first introduce some more natural notation for studying these sets.
\vskip 0.5in
\end{section}
\begin{section}{Assignments}

\vskip 0.2in
Recall the definition of the \emph{fundamental groupoid} of an unpointed topological space $X$.  We first define a \emph{path} to be a continuous function from $[0, 1]$ to $X$.  Two paths $\mathbf{P}$ and $\mathbf{Q}$ are said to be \emph{equivalent} if there is a homotopy between their corresponding functions which keeps the endpoints fixed (denoted $\mathbf{P} = \mathbf{Q}$).  Given two paths $\mathbf{P}$ and $\mathbf{Q}$ with corresponding functions $p$ and $q$ such that $p(1)=q(0)$, we may define a path $\mathbf{P}*\mathbf{Q}$ (shortened as $\mathbf{P}\mathbf{Q}$) which equals $p(2x)$ for $x \in [0, 1/2]$ and $q(2x-1/2)$ for $x \in [1/2, 1]$.  This gives us a natural definition of $\mathbf{P}^{-1}$ as the path defined by the function $p(1-x)$.
\vskip 0.3in
Then the fundamental groupoid of $X$, which we will denote $\pi_1(X)$, is precisely the collection of equivalence classes of paths in $X$, with $*$ as the group operation.  If $A \subseteq X$, we may define $\pi_1(X, A)$ as the subgroupoid consisting of paths with endpoints in $A$.
\vskip 0.3in
{\bf Definition 3.1.} Define the graph $G_{n, k}$ as follows:
\begin{itemize}
\item The vertices of this graph are the elements of $(\mathbb{Z}/k\mathbb{Z})^n$
\item For vertices $a, b$, there is an edge from $a$ to $b$ if and only if $||a-b||_{\text{sup}}=1$.  That is, we must have that every coordinate of $a-b$ is $0$ or $\pm 1$, but $a \neq b$.
\end{itemize}
\vskip 0.2in
This graph naturally embeds into the $n$-dimensional torus $\mathbb{R}^n/\mathbb{Z}^n$ by mapping $(a_1, \ldots, a_n) \in (\mathbb{Z}/k\mathbb{Z})^n$ to $(a_1/k, \ldots, a_n/k)$ and mapping the edges to geodesics.
\vskip 0.2in
Consider functions $\chi$ from paths in $G_{n, k}$ to $\mathbb{Z}^n$, with the property that if $\mathbf{P}$ and $\mathbf{Q}$ are two paths in $G_{n, k}$ which can be composed, then $\chi(\mathbf{PQ})=\chi(\mathbf{P})+\chi(\mathbf{Q})$.  Call such functions $\chi$ \emph{assignments}.  If an assignment equals zero on all contractible loops, call it \emph{closed}, and if an assignment equals zero an all loops, call it \emph{exact}.  Closed assignments on $G_{n, k}$ can be thought of as elements of $\text{Hom}(\pi_1(G_{n, k}, V), \mathbb{Z}^n)$ where $V$ is the set of vertices of $G_{n, k}$, and $G_{n, k}$ has the obvious subspace topology from $\mathbb{R}^n/\mathbb{Z}^n$.  This is because under a closed assignment, the value of a path is homotopy invariant.  Similarly, exact assignments are elements of $\text{Hom}(\pi_1(G_{n, k}, V), \mathbb{Z}^n)$ which map all loops to zero.
\vskip 0.2in
Notice that the assignments have the natural structure of an abelian group (which we call $C_{1, n}$) by addition.  The closed assignments form a subgroup $Z_{1, n} \subset C_{1, n}$, and the exact assignments form a subgroup $B_{1, n} \subset Z_{1, n}$.  We may think of an assignment as a discrete one-form on the torus, with coefficients in $\mathbb{Z}^n$, and the closed and exact assignments are analogous to closed and exact one-forms, respectively.  The quotient $Z_{1, n}/B_{1, n}$ is, as we would expect, isomorphic to $H^1(\mathbb{R}^n/\mathbb{Z}^n; \mathbb{Z}^n)$.
\vskip 0.3in
The fundamental group of $\mathbb{R}^n/\mathbb{Z}^n$ at any point is isomorphic to $\mathbb{Z}^n$.  Since it is torsion-free and abelian, we have that it is isomorphic to $H^1(\mathbb{R}^n/\mathbb{Z}^n; \mathbb{Z})$.  So, to any loop on $\mathbb{R}^n/\mathbb{Z}^n$, we can assign an element of $\mathbb{Z}^n$, called its \emph{winding vector}.  We will canonically pick this as follows.  Let $f:[0, 1] \rightarrow \mathbb{R}^n/\mathbb{Z}^n$ be such a loop, and let $\pi: \mathbb{R}^n \twoheadrightarrow \mathbb{R}^n/\mathbb{Z}^n$ be the quotient map.  Pick any point $x_0 \in \mathbb{R}^n$ such that $\pi(x_0)=f(0)$.  Then there is a unique map $\hat{f}: [0, 1] \rightarrow \mathbb{R}^n$ so that $\hat{f}(0)=x_0$ and the following diagram commutes $$\xymatrix{[0, 1]\ar@{-->}[r]^{\hat{f}}\ar[rd]_{f} & \mathbb{R}^n\ar[d]^{\pi} \\ & \mathbb{R}^n/\mathbb{Z}^n}$$ We define $w(f)=\hat{f}(0)-\hat{f}(1)$.  We call an assignment \emph{proper} if the value assigned to a loop equals its winding vector.
\vskip 0.5in
We now explain the reason for all of this notation.  Suppose that $K$ is an exactly $k$-discrete set (assume $k > 2$).  We identify the vertex $(i_1, \ldots, i_n)$ with the lattice translate of $P_m(i_1, \ldots, i_n)$, where we take $0 \le i_j \le m-1$ for $j=1, \ldots, n$.  Then we can assign a vector to each edge so that the vector corresponding to the edge going from $(i_1, \ldots, i_n)$ to $(j_1, \ldots, j_n)$ is exactly the lattice vector closest to the translation vector going from the translate of $P_m(i_1, \ldots, i_n)$ to the translate of $P_m(j_1, \ldots, j_n)$.  We can easily check that this is a proper assignment.  Similarly, any proper assignment to $G_{n, k}$ has a corresponding exactly $k$-discrete set in $\mathbb{R}^n$ (unique up to translation by $\mathbb{Z}^n$).
\vskip 0.2in
So it now follows that a subset $A \subset \mathbb{Z}^n$ is achievable if and only if for some $K$ there is an proper assignment of $G_{n, k}$ so that the value on any edge is in $A$.  We can think of the proper assignments as a coset of the group of exact assignments in the group of closed assignments.
\vskip 0.6in
\end{section}
\begin{section}{A Theorem on Connectedness}
\vskip 0.2in
We now give one more definition which allows us to state our theorem.
\vskip 0.2in
{\bf Definition 4.1.} Let $A=\{0, \pm a_1, \ldots, \pm a_k\} \subset \mathbb{Z}^n$ with $a_1, \ldots, a_k$ distinct and nonzero.  Define $\mathcal{G}(A)$, the \emph{characteristic graph} of $A$ as follows.  $\mathcal{G}(A)$ has vertices $V_1, \ldots, V_k$, and for all $i, j$, $V_i$ is connected to $V_j$ if and only if either $a_i+a_j$ or $a_i-a_j$ is in $A$.
\vskip 0.2in
{\bf Theorem 4.2.} \emph{Let $A \subset \mathbb{Z}^n$ be achievable.  Then there is some  achievable set $B \subseteq A$ such that $\mathcal{G}(B)$ is a connected component of $\mathcal{G}(A)$.}
\vskip 0.2in
We first give a formal definition of the `interior' of a subset of the torus.  We then prove a lemma which will be useful.
\vskip 0.2in
{\bf Definition 4.2.1.} Consider a compact, connected subset $Z$ of $\mathbb{R}^n/\mathbb{Z}^n$.  Let $U$ be a connected component of $\mathbb{R}^n/\mathbb{Z}^n-Z$, and let $\overline{U}$ denote the closure of $U$.  We call $U$ an \emph{interior region} of $Z$ if, for every path in $\pi_1(\overline{U}, Z)$, there is an equivalent path in $\pi_1(Z)$.  Otherwise, we call $U$ an \emph{exterior region}.
\vskip 0.3in
{\bf Lemma 4.2.2.} \emph{Consider a finite simplicial complex for $\mathbb{R}^n/\mathbb{Z}^n$.  Let $Z$ be a closed subset of $\mathbb{R}^n/\mathbb{Z}^n$ which is a union of these simplices, such that any loop in $Z$ has winding vector zero.  Then $Z$ has exactly one exterior region.}
\vskip 0.2in
{\bf Proof:} Consider the quotient map $\pi: \mathbb{R}^n \twoheadrightarrow \mathbb{R}^n/\mathbb{Z}^n$.  Let $\widehat{Z}$ denote $\pi^{-1}(Z)$.  Then $\widehat{Z}$ consists of a union of connected components, each of which is bounded (because $Z$ contains no loops of nonzero winding vector, and thus each of these components is a finite union of simplices).  We have precisely one connected component for each element of $\mathbb{Z}^n$.  Consider one such component $Z_0$.  Since $Z_0$ is bounded, $\mathbb{R}^n-Z_0$ has one unbounded region and a finite number of bounded regions.  These bounded regions clearly correspond to interior regions, because $Z_0$ is connected.  Let $U_0$ denote the unbounded component of $\mathbb{R}^n-Z_0$.
\vskip 0.3in
For each $v \in \mathbb{Z}^n$, let $Z_v$ denote the connected component of $\widehat{Z}$ which is $Z_0$, translated by $v$.  Let $U_v$ denote $U_0$'s corresponding translate.  Then $\pi(\bigcap\limits_{v \in U_v})$ is precisely the `exterior' component of $\mathbb{R}^n/\mathbb{Z}^n-Z$ we wish to prove is connected (as every other point on $\mathbb{R}^n/\mathbb{Z}^n-Z$ is in an interior region).  We then wish to show that the intersection $U=\bigcap\limits_{v \in \mathbb{Z}^n}U_v$ is connected.
\vskip 0.3in
Consider $\partial U_v \subseteq Z_v$, the boundary of $U_v$, for some $v$.  For every point $p$ on $\partial U_v$, consider a closed ball of radius $\epsilon$ around $p$, and take its intersection with $U_v$.  Now consider the union of these balls over all $p \in \partial U_v$.  Essentially, we have formed a thin `bubble' around the exterior of $Z_v$.  Label this bubble $B_{v, \epsilon}$, and let $E_{v, \epsilon}$ denote the outside boundary of this bubble: i.e., we let $E_{v, \epsilon}=\partial B_{v, \epsilon}-Z_{v}$.  If we pick $\epsilon$ to be sufficiently small, then it is clear that these $B_{v, \epsilon}$ are disjoint (over two different values $v$).
\vskip 0.3in
We claim that $B_{v, \epsilon}$ is connected for any $v$ and $\epsilon$.  It suffices to show $\partial U_v$ is connected.  We know that any point in $Z_v$ but not in $\partial U_v$ must be in the interior of some connected component of $\partial U_v$.  Additionally, it is clear that no connected component of $\partial U_v$ is in the interior of another.  Therefore, if $\partial U_v$ has two or more connected components, then there can be no path in $Z_v$ which connects them.  Hence, $\partial U_v$ is indeed connected.
\vskip 0.4in
Now suppose that we have any two points $p$ and $q$ in $U$.  If we pick $\epsilon$ sufficiently small, we may assume neither lies in any of the $B_{v, \epsilon}$.  We will show that there is a path between them in $U$.  We construct such a path according to the following algorithm.  Starting at $p$, we follow a straight line towards $q$ until we first encounter some $B_{v, \epsilon}$.  It is clear that among all points in $B_{v, \epsilon}$, the point (or points) which is closest to $q$ is in $E_{v, \epsilon}$ (see this by noting that if $r$ is such a point, then the line segment from $r$ to $q$ does not intersect $B_{v, \epsilon}$ except at $r$).  Hence, we follow a path through $B_{v, \epsilon}$ to one of the points closest to $q$.  From there, we repeat the process: again follow a straight line towards $q$, and then once we encounter another bubble $B_{v, \epsilon}$, travel within it to the point closest to $q$.  It is clear that we encounter finitely many of these `bubbles' in this process, and hence this algorithm terminates eventually at $q$.  This indeed shows that $U$ is connected, as desired.  Thus, $Z$ has precisely one exterior region, as desired.  $\blacksquare$
\vskip 0.5in
{\bf Definition 4.2.3.} Let $\ell \in \{1, 2, \ldots, n\}$.  Pick some vertex $(i_n, \ldots, i_n)$ of $G_{n, k}$ and pick distinct unit vectors $e_{j_1}, \ldots, e_{j_{\ell}}$ from the standard basis, and consider the points $(i_1, \ldots, i_n)+a_1e_{j_1}+\cdots+a_{\ell}e_{j_{\ell}}$ where $a_1, \ldots, a_{\ell} \in \{0, 1\}$.  The convex hull of these points form an $\ell$-dimensional hypercube.  We call the interior of such a hypercube an {\bf $\ell$-cell} of $G_{n, k}$.  We define the vertices of $G_{n, k}$ to be $0$-cells.
\vskip 0.2in
For an assignment $\chi$ of $G_{n, k}$, we define the \emph{support} of $\chi$ to be the union of all $\ell$-cells of $G_{n, k}$ which have an edge whose value under $\chi$ is nonzero.  We are now ready to prove Theorem 4.2.
\vskip 0.4in
\emph{Proof of Theorem 4.2.} Let $A$ be achieved by a proper assignment $\chi$ of $G_{n, k}$.  Assign a color to each connected component of $\mathcal{G}(A)$.  Correspondingly, we may color each edge of $G_{n, k}$ its corresponding color based on its value under $\chi$ (or uncolored, for edges which are assigned zero).  For each color $c \in \mathbf{C}$ (the set of colors), let $\chi_c$ be the assignment of $G_{n, k}$ which is equal to $\chi$ on edges of color $c$, and which equals zero on other edges.
\vskip 0.3in
We first claim that these $\chi_c$ have disjoint support.  Suppose otherwise, for a contradiction.  Then this implies that there is some $n$-cell which has two edges of different colors (because the support of an assignment is open).  Let the first edge have vertices $V$ and $W$ and the second have vertices $X$ and $Y$.  First suppose that $\chi(WX) \neq 0$.  Then either $\chi(VW)$ and $\chi(WX)$ correspond to connected elements in $\mathcal{G}(A)$ (because $\chi(VX)$ is their sum), or one of them is zero.  Similarly, $\chi(WX)$ and $\chi(XY)$ correspond to connected elements in $\mathcal{G}(A)$.  But this is impossible, because it would imply $\chi(VW)$ and $\chi(XY)$ are in the same connected component of $\mathcal{G}(A)$.  Hence, $\chi(WX)=0$.  Similarly, $\chi(YV)=0$.  But $\chi(VW)+\chi(WX)+\chi(XY)+\chi(YV)=0$.  Hence, $\chi(VW)+\chi(XY)=0$.  Since these two are in different connected components of $\mathcal{G}(A)$, we again have a contradiction!  It follows that the $\chi_c$ have disjoint supports.
\vskip 0.4in
Next, we claim that for all but one $c \in \mathbf{C}$, $\chi_c$ is an exact assignment.  Let $Z$ be the set of cells which are not colored.  It is clear that any loop consisting of uncolored edges is contractible.  Hence, it is clear that any loop in $Z$ has winding vector zero.  So by Lemma 4.2.2, $Z$ has exactly one exterior region.  It follows that there is a unique color $u \in \mathbf{C}$ such that the support of $\chi_u$ contains the unbounded region of $Z$, and for all other colors $c$, $\chi_c$'s support is a union of interior regions of $Z$.  We claim that these $\chi_c$ are exact assignments.
\vskip 0.25in
We prove this directly.  Pick a color $c \neq u$, and consider the assignment $\chi_c$.  Pick any loop $\mathbf{P}$ on $G_{n, k}$.  Let $C$ denote the support of $\chi_c$.  Note that $\partial C \subseteq Z$ (where, $\partial C=\overline{C}-C$).  Suppose the path starts at a point $p_0$.  We may assume that the loop begins and ends at the same point on $\partial C$ by picking a point $p_0' \in \partial C$, a path $\mathbf{Q}$ from $p_0'$ to $p_0$, and considering the path $\mathbf{Q}\mathbf{P}\mathbf{Q}^{-1}$ instead.  Thus, assume $p_0 \in \partial C$.
\vskip 0.3in
Observe that $\mathbf{P}$ corresponds to a function $P$ from $[0, 1]$ to $\mathbb{R}^n/\mathbb{Z}^n$.  Let $U$ be the set of $x \in [0, 1]$ such that $P(x)$ is in the interior of some connected component of $\partial C$.  Let $U_1, \ldots, U_m$ be the connected components of $U$.  Because $P(0), P(1) \in \partial C$, notice that these $U_i$ are open intervals.  Let us consider one such $U_i$.  Let it have endpoints $q_1$ and $q_2$, and let $\mathbf{P}_i$ be the path with associated function $P_i$ such that $P_i(x)=P(xq_2+(1-x)q_1)$.  Let $Q_1$ and $Q_2$ be the respective connected components of $\partial C$ containing $P(q_1)$ and $P(q_2)$.  It is clear that the portion of the path right after $P(q_1)$ lies in the interior of $Q_1$.  However, since $P(q_1)$ is not in the interior of any connected component of $\partial C$, it follows that $Q_1$ does not lie in the interior of any connected component of $\partial C$.  Therefore, if $q_3$ is the first point after $q_1$ such that $P(q_3) \in Q_1$, then the open interval $(q_1, q_3)$ lies in $U_i$, and $q_3 \notin U_i$.  Thus, $q_3=q_2$, and so $q_2 \in Q_1$ and $Q_1=Q_2$.  Denote this connected component of $\partial C$ by $Q$.
\vskip 0.4in
So then there is a path $\mathbf{Q}_i \in \pi_1(Q)$ such that $\mathbf{P}_i\mathbf{Q}_i^{-1} = 1$, the identity path.  Then, under $\chi_c$, $\mathbf{P}_i$ and $\mathbf{Q}_i$ evaluate to the same thing, zero.  We may, on each $U_i$, replace $\mathbf{P}|_{U_i}$ by an appropriately scaled version of $\mathbf{Q}_i$.  The result will be a path $\mathbf{Q}$ such that $\mathbf{P}\mathbf{Q}^{-1}=1$.  Moreover, no point of $\mathbf{Q}$ lies in the interior of any connected component of $\partial C$.  Since every point of $C$ lies in the interior of some connected component of $\partial C$, it follows that every edge of $\mathbf{Q}$ evaluates to $0$ under $\chi_c$.  Thus, $\mathbf{Q}$ evaluates to $0$ under $\chi_c$.  Hence, $\mathbf{P}$ does as well.
\vskip 0.4in
We now conclude the proof of the theorem as follows.  Clearly $\chi=\sum\limits_{c \in \mathbf{C}}\chi_c$.  Because the sum of a proper assignment and an exact assignment is proper, we have that $\chi_u$ is proper.  Hence, if we let $B \subseteq A$ so that $\mathcal{G}(B)$ is the connected component of $\mathcal{G}(A)$ with color $u$, then $B$ is achievable.  This proves the theorem.  $\blacksquare$
\end{section}

\vskip 0.7in
\begin{section}{A Stronger Statement and an example}
We first look at an example of a class of achievable sets.
\vskip 0.2in
{\bf Example:} \emph{Let $a_1, \ldots, a_m \in \mathbb{Z}^2$ generate $(1, 0)$ and let $b_1, \ldots, b_n \in \mathbb{Z}^2$ generate $(0, 1)$.  Consider a finite, symmetric set $A$ containing $0$, each $a_i$ and each $b_j$, and containing one of $a_i \pm b_j$ for each pair $(i, j)$.  Then $A$ is achievable.}
\vskip 0.1in
{\bf Proof.} The construction is straightforward.  Assume $\sum\limits_{i=1}^{m}a_i=(1, 0)$ and $\sum\limits_{j=1}^{n}b_j=(0, 1)$.  Consider the standard embedding of $G_{2, k}$ into $\mathbb{R}^2$.  Suppose $k$ is greater than $2m$ and $2n$.  For each $i \le m$, construct a vertical line $V_i$ which does not go through any vertices of the graph, and for each $j \le n$, construct a horizontal line $H_j$ which does not go through any vertices of the graph.  For any edge of the graph, we assign a value of the form $\sum\limits_{i=1}^{m}\delta_ia_i+\sum\limits_{j=1}^{n}\epsilon_jb_j$.  Here $\delta_i$ is equal to $0$ if the edge does not cross $V_i$, $1$ if it crosses $V_i$ from left to right, and $-1$ if it crosses $V_i$ from right to left.  Similarly, $\epsilon_j$ is $0$ if the edge does not cross $H_j$, $1$ if it crosses $H_j$ from bottom to top, and $-1$ if it crosses $H_j$ from top to bottom.  This is clearly a proper assignment.
\vskip 0.3in
Now at each intersection $P_{i, j}=V_i \cap H_j$, we perform a change depending on whether $a_i+b_j \in A$ or $a_i-b_j \in A$.  At each such intersection, there are four closest vertices: top right, top left, bottom right, and bottom left.  If we want $a_i+b_j$ in the set and not $a_i-b_j$, then we deform $V_i$ to include the top left vertex on its right side, and if we want $a_i-b_j$ in the set and not $a_i+b_j$, then we deform $V_i$ to include the top right vertex on its left side.  It is clear that once we have performed the corresponding deformation at each intersection, we will have the desired set. $\Box$
\vskip 0.3in
We can, in fact, replace $(1, 0)$ and $(0, 1)$ above by two arbitrary lattice points which generate $\mathbb{Z}^2$, by performing a corresponding affine transformation of determinant $\pm 1$.  It may be possible to do a similar construction in higher dimensions.
\vskip 0.3in
In view of  Theorem 4.2 we may ask: \emph{If we have a set $A \subset \mathbb{Z}^k$ whose characteristic graph is connected and such that $\langle A \rangle = \mathbb{Z}^k$, then is $A$ achievable?}  The answer is, in fact, no.  We will prove a stronger version of Theorem 4.2 in two dimensions, and can easily check that this implies $\{0, \pm a, \pm 2a, \pm (2a+b), \pm b\}$ is not achievable, where $\langle a, b \rangle = \mathbb{Z}^2$.
\vskip 0.4in
{\bf Lemma 5.1.} \emph{Let $k$ be a positive integer.  Let $f: [0, 1] \rightarrow \mathbb{R}^2$ be a continuous, injective function such that $f(0)=(0, 0)$ and $f(1)=(1, 0)$, and such that for any $t \in [0, 1]$, one of the two coordinates of $f(t)$ is in $(1/k)\mathbb{Z}$.  Then clearly $f$ can be extended to a unique function $\widehat{f}: \mathbb{R} \rightarrow \mathbb{R}^2$ such that $\widehat{f}$ on the interval $[n, n+1]$ is a translate by $(n, 0)$ of $f$.  Suppose that $\widehat{f}$ does not cross itself (i.e., it is injective).  Pick $d \in \mathbb{Z}$, and let $\widehat{f_1}(x)=\widehat{f}(x)+(1/d, 0)$.  Then we claim $\widehat{f}$ and $\widehat{f_1}$ take on a common value.}
\vskip 0.2in
{\bf Proof.} Because $\widehat{f}$ does not cross itself, it divides the plane into two regions.  One of these regions contains all points of the form $(x, y)$ for $y \gg 0$ (above $\widehat{f}$) and one contains all points of the form $(x, y)$ for $y \ll 0$ (below $\widehat{f}$).  Since $\widehat{f_1}$ is in one of these, without loss of generality, it is in the first.  If we define $\widehat{f_{i+1}}=\widehat{f_i}+(1/d, 0)$ for each $i \in \mathbb{N}$, we have that $\widehat{f_{i+1}}$ is above $\widehat{f_i}$ for each $i$.  The region above $\widehat{f_i}$ is clearly contained in the region above $\widehat{f}$ for each $i$.  Hence, $\widehat{f_d}$ is above $\widehat{f}$.  However, $\widehat{f}$ and $\widehat{f_d}$ have the same image.  This is a contradiction.  Thus, $\widehat{f_1}$ shares a point with $\widehat{f}$. $\Box$
\vskip 0.4in
{\bf Corollary 5.2:} \emph{Call an element of $\mathbb{Z}^n$ {\bf primitive} if its coordinates have no common divisor.  Call a loop on $G_{n, k}$ {\bf primitive} if it does not contain any nontrivial subloops.  Consider any loop $\mathbf{P} \in \pi_1(G_{n, k})$.  If its winding vector is nonzero and nonprimitive, then $\mathbf{P}$ is nonprimitive.}
\vskip 0.2in
{\bf Proof.} Assume that $\mathbf{P}$ is primitive, and suppose $\mathbf{P}$ has winding vector $(ad, bd)$ for some $d \neq 1$.  Consider the preimage $\widehat{\mathbf{P}}$ in $\mathbb{R}^2$.  Without loss of generality, $\mathbf{P}$ passes through the origin on the torus, and hence, this $\widehat{\mathbf{P}}$ passes through every lattice point.  One connected component of $\widehat{\mathbf{P}}$ passes through every multiple of $(a, b)d$, and another passes through every point that is the sum of $(a, b)$ and a multiple of $(a, b)d$.  By assumption, neither of these components crosses itself.  By the above lemma (after an affine transformation), these two components must intersect, which means $\mathbf{P}$ crosses itself somewhere between its two endpoints. $\blacksquare$
\vskip 0.5in
{\bf Theorem 5.3.} \emph{Let $A \subset \mathbb{Z}^2$ be achievable.  Let $S_1, \Delta_1, \ldots, S_m, \Delta_m \subset A$ be symmetric sets such that}
\begin{itemize}
\item \emph{$\bigcup\limits_{i=1}^{m}(S_i \cup \Delta_i)=A-\{0\}$}

\item \emph{$\Delta_i$ is the set of elements of $A$ which are adjacent to elements of $S_i$ in $\mathcal{G}(A)$, but which are not in $S_i$}

\item \emph{$\langle \Delta_i \rangle$ does not contains any primitive elements}

\item \emph{The $S_i$ are disjoint, and no two distinct $S_i$ and $S_j$ contains respective elements $a_i$ and $a_j$ which are adjacent in $\mathcal{G}(A)$}

\end{itemize}
\emph{Then either $S_i \cup \langle \Delta_i\rangle$ is achievable for some $i$, or $\bigcup\limits_{i=1}^{m}\langle \Delta_i \rangle$ is achievable.}
\vskip 0.1in
{\bf Proof.} The method of proof is similar to that in Theorem 4.2.  If $A$ is achievable, then there is a proper assignment $\chi$ of $G_{n, k}$ which achieves it.  Similar to before, we let $\chi_i$ denote the assignment which equals $\chi$ on edges where $\chi$ takes value in $S_i$, and equals zero on all other edges.  Of course in this case, $\chi_i$ is not necessarily exact or proper, but is, as we will see, `almost' exact or `almost' proper, so we can still use an approach similar to our previous one.
\vskip 0.4in
Let $C_i$ denote the boundary of the support of $\chi_i$ (i.e., $\overline{\text{Supp}(\chi_i)}-\text{Supp}(\chi_i)$).  Consider any edge in $C_i$.  We can easily check that its value under $\chi$ is either an element of $\Delta_i$, or is a sum of two elements of $\Delta_i$.  Hence, it follows that if $\mathbf{P}$ is a path in $C_i$, then $w(\mathbf{P}) \in \langle \Delta_i \rangle$.  However, $\langle \Delta_i \rangle$ contains no primitive elements, and so it follows that $C$ contains no loops with nonzero winding vector.  Lemma 4.2.2 then tells us that $C$ has exactly one exterior region.  We consider two cases: the first being that $\text{Supp}(\chi_i)$ is the exterior region of $C_i$ for some $i$, and the second being that it is an interior region for every $i$.
\vskip 0.4in
Suppose that $\text{Supp}(\chi_i)$ contains the exterior region of $C_i$ for some $i$.  We then claim that $S_i \cup \langle \Delta_i \rangle$ is achievable.  Consider all edges of $G_{n, k}$ which are in some interior region of $C_i$, and let $\psi$ be the assignment which equals $\chi$ on these edges, and zero on all other edges.
\vskip 0.3in
Recall that an assignment of $G_{n, k}$ is a function from paths to elements of $\mathbb{Z}^n$ so that composing two paths adds their values.  By composing such a function with the quotient map $\mathbb{Z}^n \rightarrow \mathbb{Z}^n/\langle \Delta_i \rangle$, we can obtain a function from paths to elements of $\mathbb{Z}^n/\langle \Delta_i \rangle$.  We call this a \emph{$\Delta_i$-assignment}.  As before, we define \emph{closed} $\Delta_i$-assignments to be elements of $\text{Hom}(\pi_1(G_{n, k}, A), \mathbb{Z}^n/\langle \Delta_i \rangle)$, and we can define \emph{exact} $\Delta_i$-assignments to be closed and zero on all loops.
\vskip 0.3in
Now by composing the map $\mathbb{Z}^n \rightarrow \mathbb{Z}^n/\langle \Delta_i \rangle$ with $\psi$, we get a $\Delta_i$-assignment $\tilde{\psi}$.  We claim that this is exact.  Notice that because $\chi$ takes on values in $\langle \Delta_i \rangle$ for all edges on $C_i$, this claim follows from an argument identical to that in Theorem 4.2.  We now show that we can change $\psi$ slightly to an assignment $\psi'$ which is exact.
\vskip 0.3in
We will change the value of $\psi$ only on edges in the interior of $C_i$, and only then by elements of $\langle \Delta_i \rangle$.  Notice that $\overline{\text{Supp}(\psi)}$ contains no loops with nonzero winding vector.  Thus, if we have any two vertices in the same connected component of $\overline{\text{Supp}(\psi)}$, and consider two paths between them, both contained in $\text{Supp}(\psi)$, then these two paths are homotopic in $\text{Supp}(\psi)$, and thus take the same value under $\psi$.
\vskip 0.3in
Pick a vertex on $P$ on the boundary of $\text{Supp}(\psi)$.  Consider any other vertex $Q$ on the same connected component of this boundary.  Then any two paths connecting them in $\overline{\text{Supp}(\psi)}$ take the same value in $\chi$.  Hence, we may assign each such point $Q$ the corresponding value $f(Q)$ (where $f(P)=0$).  We then modify $\psi$ as follows.  For every pair of vertices $R \in \text{Supp}(\psi)$ and $Q \in \partial \text{Supp}(\psi)$ with an edge from $Q$ to $R$, we increase the value of this edge by $f(Q)$.  We perform this change for all such pairs $Q$ and $R$.
\vskip 0.4in
It is clear that the resulting assignment $\psi'$ is a lift of $\tilde{\psi}$, so it suffices to show that it is exact.  Consider any loop $\mathbf{P}$ on $G_{n, k}$.  Its intersection with $\overline{\text{Supp}(\psi)}$ is a collection of closed intervals.  We will prove that $\mathbf{P}$ takes the value $0$ by induction on the number of such closed intervals.  If $\mathbf{P}$ does not intersect $\text{Supp}(\psi)$, then it clearly takes the value $0$.  Now suppose that $\mathbf{P}$ intersects $\text{Supp}(\psi)$ in some closed interval with endpoints $P$ and $Q$.  Let $\mathbf{S}$ denote the portion of the path from $P$ to $Q$, and let $\mathbf{T}$ be a path from $P$ to $Q$ along $\partial \text{Supp}(\psi)$ (clearly $P$ and $Q$ are in the same connected component).  It is clear that $\psi'(\mathbf{S})=\psi'(\mathbf{T})$, because $\psi'(\mathbf{S}\mathbf{T}^{-1})=\psi(\mathbf{S}\mathbf{T}^{-1})=0$.  Hence, we may replace this portion $\mathbf{S}$ of $\mathbf{P}$ with $\mathbf{T}$, and the resulting path, under $\psi'$, takes the value $\psi'(\mathbf{P})$.  This new path intersects $\overline{\text{Supp}(\psi)}$ in one fewer closed interval, so this completes the inductive step, and thus shows $\psi'$ is exact.  Thus, $\chi-\psi'$ is proper and every edge takes value in $S_i \cup \langle \Delta_i \rangle$.

\vskip 0.5in
Now let us consider the case where $\text{Supp}(\chi_i)$ is contained in the interior regions of $C_i$ for every $i$.  As above, we may let $\phi_i$ denote the $\Delta_i$-assignment generated by $\chi_i$ for each $i$, and we notice that each $\phi_i$ is exact.  We then, in a similar manner, construct assignments $\chi_i'$ which are exact and only slightly different from the $\chi_i$.  We then may similarly consider the proper assignment $\chi-\sum\limits_{i}^{m}\chi_i'$ whose edge labels are all in $\bigcup\limits_{i=1}^{m}\langle \Delta_i \rangle$. $\blacksquare$
\vskip 0.5in
We conjecture the following generalization of Theorem 5.3.
\vskip 0.3in
{\bf Conjecture 5.4.} \emph{Let $A \subset \mathbb{Z}^n$ be achievable.  Let $S_1, \Delta_1, \ldots, S_m, \Delta_m \subset A$ be symmetric sets such that}
\begin{itemize}
\item \emph{$\bigcup\limits_{i=1}^{m}(S_i \cup \Delta_i)=A-\{0\}$}

\item \emph{$\Delta_i$ is the set of elements of $A$ which are adjacent to elements of $S_i$ in $\mathcal{G}(A)$, but which are not in $S_i$}

\item \emph{$\mathbb{Z}^n/\langle \Delta_i \rangle$ is not cyclic}

\item \emph{The $S_i$ are disjoint, and no two distinct $S_i$ and $S_j$ contain respective elements $a_i$ and $a_j$ which are adjacent in $\mathcal{G}(A)$}

\end{itemize}
\emph{Then either $S_i \cup \langle \Delta_i\rangle$ is achievable for some $i$, or $\bigcup\limits_{i=1}^{m}\langle \Delta_i \rangle$ is achievable.}

\end{section}

\vskip 0.7in
\begin{section}{A More General Question and a Conjecture}
\vskip 0.3in
\hskip 0.23in {\bf Definition 6.1.} Let $G$ be a group.  Consider two finite subsets $S_1$ and $S_2$.  Let us say that $S_1 \le S_2$ if $S_1g \subseteq S_2$ for some $g \in G$, and $S_1=S_2$ if $S_1g=S_2$.  We can easily check that this property is transitive, and thus defines a partial order on finite subsets of $G$, modulo this equivalence relation.  Call a poset $\mathcal{A}$ of such finite subsets of $G$ \emph{ideal} if $\mathcal{A}$ consists of equivalence classes of sets and for each $S$ in $\mathcal{A}$, $S' \in \mathcal{A}$ for each $S' \le S$.  It is clear that any ideal is graded by the sizes of these subsets.  If such a set $S \subseteq G$ has $n+1$ elements, we say it has \emph{rank} $n$.
\vskip 0.4in
Now suppose that $G$ acts geometrically on a metric space $X$.  Call an ideal $\mathcal{A}$ \emph{achievable} if $\exists K \subset X$ such that
\vskip 0.1in
$\bullet$ $K$ is compact and $G(K)=X$ (i.e., $K$ is an $\mathcal{N}$-\emph{set}).
\vskip 0.1in
$\bullet$ For any $p \in X$, there is some $S \in \mathcal{A}$ such that $S(p)=G(p) \cap K$.
\vskip 0.3in
In general, suppose that $K \subseteq X$ is an $\mathcal{N}$-\emph{set}, and $\mathcal{A}$ is an ideal.  We say that $K$ \emph{achieves} $\mathcal{A}$ if, for any $p \in X$, the set $\{g: g(p) \in C\}$ is in $\mathcal{A}$.  We denote the achieved poset of $C$ by $\mathcal{A}(K)$.
\vskip 0.4in
{\bf Lemma 6.2.1.} \emph{Let $K_0 \supseteq K_1 \supseteq \cdots \supseteq K$ be a sequence of compact sets in $X$ such that $\bigcap\limits_{i=0}^{\infty}K_i=K$.  Then for some $m$, $\mathcal{A}(K_m)=\mathcal{A}(K)$.}
\vskip 0.2in
{\bf Proof.} Suppose otherwise.  Let $\mathcal{A}_i$ denote $\mathcal{A}(K_i)$ for $i=0, 1, \ldots$.  By an argument similar to that in Lemma 2.1, there is some $S \subseteq G$ which is contained in every $\mathcal{A}_i$, but not contained in $\mathcal{A}(K)$.  Hence, for $i=0, 1, \ldots$, we have $p_i \in K_i$ so that $S(p_i) \subseteq K_i$.
\vskip 0.2in
Let $S$ have rank $r$.  By an argument similar to that in Lemma 2.1, the sequence $\{Sp_i\}_{i \ge 0}$ has a convergent subsequence $Sp_{a_1}, Sp_{a_2}, \ldots$.  Call its limit $\{q_0, \ldots, q_r\}$.
\vskip 0.3in
Let $g_1, \ldots, g_r$ be the elements of $S$.  Because this sequence converges under the supremum norm, we can take $q_0$ to be the limit of $p_{a_1}, p_{a_2}, \ldots$, and $q_i$ to be the limit of $g_ip_{a_1}, g_ip_{a_2}, \ldots$ for each $i$.  Because each of $q_0, \ldots, q_r$ is in $K_0, \ldots$, it is in $K$.  By an argument similar to that in Lemma 2.1, $g_iq_0=q_i$, and hence, $S \in \mathcal{A}(K)$, a contradiction.  Thus, no such $S$ exists, and so $\mathcal{A}_m=\mathcal{A}(K)$ for some $m$ sufficiently large.

\vskip 0.5in
Now let us suppose that $X=\mathbb{R}^n$ and $G=\mathbb{Z}^n$.
\vskip 0.2in
{\bf Corollary 6.2.2.} \emph{If $\mathcal{A}$ is achievable, then it is achievable by an exactly discrete set.}
\vskip 0.1in
{\bf Proof.} Again by Lemma 6.2.1, it is clear that $\mathcal{A}$ is achievable by a discrete set.
\vskip 0.2in
Now we must show that $\mathcal{A}$ is achievable by an exactly discrete set.  First, it is clear that some subset of $\mathcal{A}$ is achievable by an exactly discrete set.  Suppose now that this subset is $\mathcal{Q}$, and we want to achieve $\mathcal{Q} \cup \Lambda_S$ for some $S \subset \mathbb{Z}^n$, where $\Lambda_S$ is the smallest ideal containing $S$.  Let $S=(s_0, s_1, \ldots, s_k)$, where $k<2^n$.  Pick one hypercube cell in our discrete set, and cut it up into $4^n$ smaller cubes.  Look at the inner $2^n$ of these, and for each $i=1, 2, \ldots, k$, translate one by $s_i-s_0$.  Since every additional set of lattice points which is achieved by this modification is a subset of $S$, and $S$ is achieved in this way, it follows that this new modified exactly discrete set achieves $\mathcal{Q} \cup \Lambda_S$.  A sequence of such modifications will give us $\mathcal{A}$.  $\blacksquare$
\vskip 0.5in
{\bf Conjecture 6.3:} \emph{Suppose that $\mathcal{A}$ is an achievable poset of sets of points in $\mathbb{Z}^n$.  Then there is an achievable poset $\mathcal{A}' \subseteq \mathcal{A}$ such that
\vskip 0.1in
$\bullet$ The maximal rank in $\mathcal{A}'$ is $n$, and every element has a rank $n$ ancestor
\vskip 0.1in
$\bullet$ Call two rank $n-1$ elements {\bf connected} if they have a common rank $n$ parent.  Then the rank $n-1$ elements of $\mathcal{A}'$ are in one connected component.}
\vskip 0.3in
This conjecture is a natural generalization of Theorem 4.2.  It still is not clear if this conjecture could be modified to an even stronger analogue of Theorem 5.3.  We hope that the general ideas in the proof of Theorem 4.2 might carry over.  Below we sketch a general idea of how the proof might proceed.  We can first prove that if $\mathcal{A}$ is achievable, then the poset consisting of all elements of $\mathcal{A}$ of rank at most $n$ is achievable (we do so essentially using the topological definition of dimension).  Now assign a color to each connected component of the set of rank $n-1$ elements of $\mathcal{A}$.  We first show the following statement.
\vskip 0.2in
{\bf Proposition:} \emph{Let $\mathcal{A}$ be an achievable poset.  Then $\mathcal{A}$ is achieved by a proper assignment of $G_{n, k}$, for some $k$, such that for every $n$-cell containing $n$ level sets in this assignment, some $(n-1)$-cell face of this $n$-cell contains these $n$ level sets.}
\vskip 0.2in
Hence, each $(n-1)$-cell inherits a unique color.  Call two $(n-1)$-cells \emph{strongly connected} if they are contains within a common $n$-cell.  Then we notice that if two $(n-1)$-cells are strongly connected, they must either be the same color, or one must be uncolored.  If we let $U$ be the set of uncolored $(n-1)$-cells, we have that the other $(n-1)$-cells are split into strongly connected components which are monochromatic.
\vskip 0.2in
We now wish to show that we can reduce to the case where there is only one strongly connected component.  We need to show that we can extend the assignment on the edges of $U$ to all but one connected component of $\mathbb{R}^n/\mathbb{Z}^n-U$ such that no new elements are added to $\mathcal{A}$ in the process.  We do not know how to approach this step yet, or if it is possible as a method of proof.
\end{section}
\vskip 1in

\begin{section}{Conclusion and Acknowledgements}
\vskip 0.1in
The theorems we have proven still do not give a complete characterization of achievable sets even in $\mathbb{Z}^2$, so we may ask if such a characterization exists.  We also wish for proofs of Conjectures 5.4 and 6.2.
\vskip 0.2in
We can state Nathanson's question more generally for geometric group actions as well.  Suppose that $X$ is a metric space, and $G$ is a group which acts geometrically on $X$.  That is, $X$ must be boundedly compact and geodesic, and the group action must be properly discontinuous, co-compact, and isometric.  Call a subset $K$ of $X$ an $\mathcal{N}$-\emph{set} if $K$ is compact and covers every orbit of $G$.  Then define $$A(K):=\{a \in G: K \cap aK \neq \emptyset\}$$ to be the \emph{achieved set} of $K$.  Conversely, call a set $A \subseteq G$ \emph{achievable} if $A=A(K)$ for some $\mathcal{N}$-\emph{set} $K$.  Then what sets are achievable?  As in Section 6, we can also define achievable posets, and ask what ideals of $\mathcal{P}(G)$ are achievable.  Perhaps we could generalize the theorems above to certain classes of manifolds and groups.
\vskip 0.2in
This research was done over the summer at the University of Minnesota Duluth, with grants from the NSF (grant number DMS 0754106) and the NSA (grant number H98230-09-1-0115).  I would like to thank Joe Gallian for running the program and for his advice during my research.  I would also like to thank my advisors, Nathan Pflueger and Maria Monks for their help in my research.  I especially would like to thank Pflueger and Gallian for their help with this paper.
\vskip 0.2in
\end{section}


\begin{thebibliography}{9}
\bibitem{1} L. A. Borisov and R. Jin, \emph{Finding integral diagonal pairs in a two dimensional $\mathcal{N}$-\emph{set}}, arXiv:1007.1441v1

\bibitem{2} Z. Ljujic and C. Sanabria, \emph{The inverse problem for the lattice points}, arXiv:1007.1782v1

\bibitem{3} M. B. Nathanson, \emph{An inverse problem in number theory and geometric group theory}, arXiv:0901.1458v2

\end{thebibliography}
\end{document}